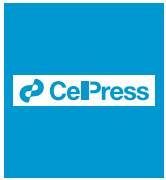

Contents lists available at ScienceDirect

# Heliyon

journal homepage: www.cell.com/heliyon

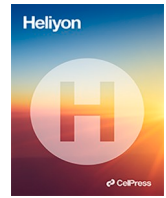

Research article

# Modified conformable double Laplace–Sumudu approach with applications

Shams A. Ahmed [a,b], Rania Saadeh [c,*], Ahmad Qazza [c], Tarig M. Elzaki [d]

[a] *Department of Mathematic, Jouf University, Tubarjal, Saudi Arabia*
[b] *Department of Mathematic, University of Gezira, Sudan*
[c] *Department of Mathematic, Zarqa University, Zarqa, 13110, Jordan*
[d] *Department of Mathematic, Alkamil, Jeddah, University of Jeddah, Kingdom of Saudi Arabia*

ABSTRACT

In this study, we combine two novel methods, the conformable double Laplace-Sumudu transform (CDLST) and the modified decomposition technique. We use the new approach called conformable double Laplace-Sumudu modified decomposition (CDLSMD) method, to solve some nonlinear fractional partial differential equations. We present the essential properties of the CDLST and produce new results. Furthermore, five interesting examples are discussed and analyzed to show the efficiency and applicability of the presented method. The results obtained show the strength of the proposed method in solving different types of problems.

## 1. Introduction

Fractional partial differential equations have recently been shown to be crucial for describing a wide range of phenomena and applications in science and engineering, including fluid dynamics, electrical circuits, optics, physics, and others [1–10]. In the literature, a number of definitions of fractional integrals and derivatives, including those by Caputo, Hadamard, and others, have been explored [11–17]. For instance, the quotient and product rules [18–23] are among the features of the ordinary derivatives that are not always met by fractional derivatives. A novel definition of a conformable fractional derivative, which virtually meets most of the standard features of derivatives, was proposed by the researchers in Ref. [24,25].

Numerous mathematicians and authors have created new methods to solve conformal differential problems, including the simplest method [26], Kudryashov method [27], double Shehu transform [28], Tanh method [29,30], reliable methods [31], double Sumudu transform [32,33], conformable Laplace transform (CLT) method [34,35], conformable double Laplace transform method [36], and conformable Sumudu transform (CST) method and others [37–44].

The double Laplace-Sumudu transformation method, a unique approach to double transformations, has been introduced recently and has been successfully used to solve families of fractional partial differential equations [45–48]. Unfortunately, this integral transformation cannot directly handle nonlinear equations, unlike other integral transformations. Because of this, mathematicians have developed new techniques that incorporate integral transformations with numerical methods like the decomposition method, the perturbation method, the iteration variation method, and others [49–59].

In this study, we introduce a unique method for numerically solving nonlinear conformable fractional differential equations, known as the CDLSMD method. The approach is effective in producing approximate series solutions that converge. With no evidence of noise terms, the approach suggested in this work quickly converges the obtained series solutions. While the components of the series solutions are represented by other analysis techniques, such as the homotopy perturbation approach and the Adomian decomposition method, these terms should also have noise terms [60].

Following that, we introduce a nonlinear conformable partial differential equation that is given by

$$b\frac{\partial^{n\gamma}}{\partial y^{n\gamma}}\left(\psi\left(\frac{x^{\eta}}{\eta},\frac{y^{\gamma}}{\gamma}\right)\right) + c\frac{\partial^{m\eta}}{\partial x^{m\eta}}\left(\psi\left(\frac{x^{\eta}}{\eta},\frac{y^{\gamma}}{\gamma}\right)\right) + \mathcal{N}\left[\psi\left(\frac{x^{\eta}}{\eta},\frac{y^{\gamma}}{\gamma}\right)\right] = g\left(\frac{x^{\eta}}{\eta},\frac{y^{\gamma}}{\gamma}\right), \frac{x^{\eta}}{\eta},\frac{y^{\gamma}}{\gamma} > 0, 0 < \gamma, \eta \leq 1, m,n \in \mathbb{N}, \quad (1)$$

with (1), we associate the initial conditions

$$\frac{\partial^{j\gamma}\psi\left(\frac{x^{\eta}}{\eta},0\right)}{\partial y^{j\gamma}} = f_j\left(\frac{x^{\eta}}{\eta}\right), j = 0, 1, ..., n-1, \quad (2)$$

and the conditions

$$\frac{\partial^{k\eta}\psi\left(0,\frac{y^{\gamma}}{\gamma}\right)}{\partial x^{k\eta}} = h_k\left(\frac{y^{\gamma}}{\gamma}\right), k = 0, 1, ..., m-1. \quad (3)$$

here, $c$ are constants, $\mathcal{N}$ is the nonlinear operator, and $g\left(\frac{x^{\eta}}{\eta},\frac{y^{\gamma}}{\gamma}\right)$ is a given function that can be expressed as $g\left(\frac{x^{\eta}}{\eta},\frac{y^{\gamma}}{\gamma}\right) = g_1\left(\frac{x^{\eta}}{\eta},\frac{y^{\gamma}}{\gamma}\right) + g_2\left(\frac{x^{\eta}}{\eta},\frac{y^{\gamma}}{\gamma}\right)$.

In this work, we introduce the CDLSMD method to investigate certain applications similar to those in (1) connected to conditions (2) and (3). The modified decomposition method [60–64] and the CDLST method are combined in the CDLSMD approach. Contrary to other analysis approaches, abandoning the CDLST approach along with the modified decomposition method enables for quick convergence to the exact solutions without discretization or linearization, as shown in Ref. [60]. This study's major objective is to propose the CDLSMD approach for investigating how conformal nonlinear partial differential equations are solved.

The research presents the fundamental definitions and facts of the CDLST, then the key concept of the novel technique, and lastly the solution to certain numerical applications to demonstrate the effectiveness and adaptability of the proposed approach.

## 2. Fractional derivatives

**Definition 1**. [24]. Given $\psi : (0,\infty) \to R$, then the CFD (conformable fractional derivative) of order $\gamma$ of $\psi$ is defined by

$$\frac{d^{\gamma}}{dy^{\gamma}}\psi(y) = \lim_{\sigma \to 0}\frac{\psi(y + \sigma y^{1-\gamma}) - \psi(y)}{\sigma}, y > 0, \sigma \in \left(0,1\right].$$

**Definition 2**. [44] Given $\psi\left(\frac{x^{\eta}}{\eta},\frac{y^{\gamma}}{\gamma}\right) : R^+ \times R^+ \to R$, then the conformable fractional partial derivatives (CFPDs) of order $\eta$ and $\gamma$ of the function $\psi\left(\frac{x^{\eta}}{\eta},\frac{y^{\gamma}}{\gamma}\right)$ is defined by

$$\frac{\partial^{\eta}}{\partial x^{\eta}}\psi\left(\frac{x^{\eta}}{\eta},\frac{y^{\gamma}}{\gamma}\right) = \lim_{\rho \to 0}\frac{\psi\left(\frac{x^{\eta}}{\eta} + \rho x^{1-\eta},\frac{y^{\gamma}}{\gamma}\right) - \psi\left(\frac{x^{\eta}}{\eta},\frac{y^{\gamma}}{\gamma}\right)}{\rho},$$

$$\frac{\partial^{\gamma}}{\partial y^{\gamma}}\left(\frac{x^{\eta}}{\eta},\frac{y^{\gamma}}{\gamma}\right) = \lim_{\sigma \to 0}\frac{\psi\left(\frac{x^{\eta}}{\eta},\frac{y^{\gamma}}{\gamma} + \sigma y^{1-\gamma}\right) - \psi\left(\frac{x^{\eta}}{\eta},\frac{y^{\gamma}}{\gamma}\right)}{\sigma},$$

where $\frac{x^{\eta}}{\eta},\frac{y^{\gamma}}{\gamma} > 0, 0 < \eta,\gamma \leq 1$, $\frac{\partial^{\eta}}{\partial x^{\eta}}$ and $\frac{\partial^{\gamma}}{\partial y^{\gamma}}$ denote the $\eta$ and $\gamma$ fractional derivatives respectively.

**Definition 3**. Let $< \gamma \leq m+1$, $m = 0,1,2,...$, and set $\alpha = \gamma - m$ then the conformable fractional integral (CFI), starting $a$ of order $\gamma$ of $\psi(y)$ is defined by

$$\left(I_{\gamma}^{a}\right)\psi(y) = \frac{1}{m!}\int_{a}^{y}(y-\tau)^{m}(\tau-a)^{\alpha-1}\psi(\tau)d\tau,$$

where $\psi(\tau)$ is continuous

The following theorems present some results related to partial derivatives the conformable partial fractional derivatives, as follows.

**Theorem 4**. [39]. *suppose* $0 < \eta,\gamma \leq 1$ *and* $\psi\left(\frac{x^{\eta}}{\eta},\frac{y^{\gamma}}{\gamma}\right)$ *be a differentiable at a point* $\frac{x^{\eta}}{\eta},\frac{y^{\gamma}}{\gamma} > 0$, *then*

$$\frac{\partial^\eta}{\partial x^\eta}\psi\left(\frac{x^\eta}{\eta},\frac{y^\gamma}{\gamma}\right) = x^{-\eta+1}\frac{\partial\psi\left(\frac{x^\eta}{\eta},\frac{y^\gamma}{\gamma}\right)}{\partial x}.$$

$$\frac{\partial^\gamma}{\partial y^\gamma}\psi\left(\frac{x^\eta}{\eta},\frac{y^\gamma}{\gamma}\right) = y^{-\gamma+1}\frac{\partial\psi\left(\frac{x^\eta}{\eta},\frac{y^\gamma}{\gamma}\right)}{\partial y}. \quad (4)$$

Proof. We find, from the definition 1,

$$\frac{\partial^\eta}{\partial x^\eta}\psi\left(\frac{x^\eta}{\eta},\frac{y^\gamma}{\gamma}\right) = \lim_{\rho\to 0}\frac{\psi\left(\frac{x^\eta}{\eta}+\rho\mathscr{X}^{1-\eta},\frac{y^\gamma}{\gamma}\right)-\psi\left(\frac{x^\eta}{\eta},\frac{y^\gamma}{\gamma}\right)}{\rho},$$

putting $\rho x^{1-\eta} = \theta$.

$$\frac{\partial^\eta}{\partial x^\eta}\psi\left(\frac{x^\eta}{\eta},\frac{y^\gamma}{\gamma}\right) = \lim_{\theta\to 0}\frac{\psi\left(\frac{x^\eta}{\eta}+\theta,\frac{y^\gamma}{\gamma}\right)-\psi\left(\frac{x^\eta}{\eta},\frac{y^\gamma}{\gamma}\right)}{\theta x^{\eta-1}} = x^{-\eta+1}\lim_{\theta\to 0}\frac{\psi\left(\frac{x^\eta}{\eta}+\theta,\frac{y^\gamma}{\gamma}\right)-\psi\left(\frac{x^\eta}{\eta},\frac{y^\gamma}{\gamma}\right)}{\theta} = x^{-\eta+1}\frac{\partial\psi\left(\frac{x^\eta}{\eta},\frac{y^\gamma}{\gamma}\right)}{\partial x}$$

By similar arguments, we prove (4).

**Theorem 5.** *Suppose $0 < \eta,\gamma \leq 1$ and $c, d, m_1, m_2, \delta$, and $\varepsilon \in \mathbb{R}$; then*

$$\frac{\partial^\eta}{\partial x^\eta}\left(c\psi\left(\frac{x^\eta}{\eta},\frac{y^\gamma}{\gamma}\right)+d\zeta\left(\frac{x^\eta}{\eta},\frac{y^\gamma}{\gamma}\right)\right) = c\frac{\partial^\eta}{\partial x^\eta}\psi\left(\frac{x^\eta}{\eta},\frac{y^\gamma}{\gamma}\right)+d\frac{\partial^\eta}{\partial x^\eta}\zeta\left(\frac{x^\eta}{\eta},\frac{y^\gamma}{\gamma}\right),$$

$$\frac{\partial^\eta}{\partial x^\eta}\left(e^{\delta\frac{x^\eta}{\eta}+\varepsilon\frac{y^\gamma}{\gamma}}\right) = \delta e^{\delta\frac{x^\eta}{\eta}+\varepsilon\frac{y^\gamma}{\gamma}},$$

$$\frac{\partial^\gamma}{\partial y^\gamma}\left(e^{\delta\frac{x^\eta}{\eta}+\varepsilon\frac{y^\gamma}{\gamma}}\right) = \varepsilon e^{\delta\frac{x^\eta}{\eta}+\varepsilon\frac{y^\gamma}{\gamma}},$$

$$\frac{\partial^\eta}{\partial x^\eta}\left(\frac{x^\eta}{\eta}\right)^{m_1}\left(\frac{y^\gamma}{\gamma}\right)^{m_2} = m_1\left(\frac{x^\eta}{\eta}\right)^{m_1-1}\left(\frac{y^\gamma}{\gamma}\right)^{m_2},$$

$$\frac{\partial^\gamma}{\partial y^\gamma}\left(\frac{x^\eta}{\eta}\right)^{m_1}\left(\frac{y^\gamma}{\gamma}\right)^{m_2} = m_2\left(\frac{x^\eta}{\eta}\right)^{m_1}\left(\frac{y^\gamma}{\gamma}\right)^{m_2-1},$$

$$\frac{\partial^\eta}{\partial x^\eta}\left(\sin\left(\frac{x^\eta}{\eta}\right)\sin\left(\frac{y^\gamma}{\gamma}\right)\right) = \cos\left(\frac{x^\eta}{\eta}\right)\sin\left(\frac{y^\gamma}{\gamma}\right),$$

$$\frac{\partial^\gamma}{\partial y^\gamma}\left(\sin\left(\frac{x^\eta}{\eta}\right)\sin\left(\frac{y^\gamma}{\gamma}\right)\right) = \sin\left(\frac{x^\eta}{\eta}\right)\cos\left(\frac{y^\gamma}{\gamma}\right).$$

The proof can be obtained using Theorem 4.

## 3. Basic facts of the CDLST

**Definition 6.** [65] Assume that $\psi\left(\frac{x^\eta}{\eta},\frac{y^\gamma}{\gamma}\right), \frac{x^\eta}{\eta},\frac{y^\gamma}{\gamma} \in \mathbb{R}^+$ and $0 < \eta,\gamma \leq 1$ is a function of two variables. Then,

(i) The fractional conformable Laplace transform of $\psi\left(\frac{x^\eta}{\eta}\right)$ order $\eta$, denoted by $\mathscr{L}_x^\eta\left[\psi\left(\frac{x^\eta}{\eta}\right)\right] = \Psi_\eta(v)$ and defined as

$$\mathscr{L}_x^\eta\left[\psi\left(\frac{x^\eta}{\eta}\right)\right] = \Psi_\eta(v) = \int_0^\infty e^{-v\frac{x^\eta}{\eta}}\psi\left(\frac{x^\eta}{\eta}\right)x^{\eta-1}dx, x > 0.$$

(ii) The conformable fractional Sumudu transform of $\psi\left(\frac{y^\gamma}{\gamma}\right)$ order $\gamma$, denoted by $\mathscr{S}_y^\gamma\left[\psi\left(\frac{y^\gamma}{\gamma}\right)\right] = \Psi_\gamma(\omega)$ and defined as

$$\mathscr{S}_y^\gamma\left[\psi\left(\frac{y^\gamma}{\gamma}\right)\right] = \Psi_\gamma(\omega) = \frac{1}{\omega}\int_0^\infty e^{-\frac{1}{\omega}\frac{y^\gamma}{\gamma}}\psi\left(\frac{y^\gamma}{\gamma}\right)\frac{y^{\gamma-1}}{\gamma}dy, y > 0.$$

(iii) The CDLST of $\psi\left(\frac{x^\eta}{\eta}, \frac{y^\gamma}{\gamma}\right)$, denoted by $\mathscr{L}_x^\eta \mathscr{S}_y^\gamma\left[\psi\left(\frac{x^\eta}{\eta}, \frac{y^\gamma}{\gamma}\right)\right] = \Psi_{\eta,\gamma}(v, \omega)$ and defined as

$$\mathscr{L}_x^\eta \mathscr{S}_y^\gamma\left[\psi\left(\frac{x^\eta}{\eta}, \frac{y^\gamma}{\gamma}\right)\right] = \Psi_{\eta,\gamma}(v, \omega) = \frac{1}{\omega} \int_0^\infty \int_0^\infty e^{-v\frac{x^\eta}{\eta} - \frac{1}{\omega}\frac{y^\gamma}{\gamma}} \psi\left(\frac{x^\eta}{\eta}, \frac{y^\gamma}{\gamma}\right) x^{\eta-1} y^{\gamma-1} dx\, dy,$$

where $v, \omega \in \mathbb{C}$ are Laplace – Sumudu variables of $\frac{y^\gamma}{\gamma}, \frac{x^\eta}{\gamma}$ and $0 < \gamma, \eta \leq 1$

The inverse CDLST $\mathscr{L}_x^{-1} \mathscr{S}_y^{-1}[\Psi_{\eta,\gamma}(v, \omega)] = \psi\left(\frac{x^\eta}{\eta}, \frac{y^\gamma}{\gamma}\right)$ is defined by

$$\mathscr{L}_x^{-1} \mathscr{S}_y^{-1}[\Psi_{\eta,\gamma}(v, \omega)] = \psi\left(\frac{x^\eta}{\eta}, \frac{y^\gamma}{\gamma}\right) = \frac{1}{2\pi i} \int_{\alpha-i\infty}^{\alpha+i\infty} e^{v\frac{x^\eta}{\eta}} \left( \frac{1}{2\pi i} \int_{\beta-i\infty}^{\beta+i\infty} \frac{1}{\omega} e^{\frac{y^\gamma}{\omega\gamma}} \Psi_{\eta,\gamma}(v, \omega) d\omega \right) dv.$$

The usual double Laplace–Sumudu transform and the CDLST are related by the following theorem.

**Theorem 7.** *Suppose that $c, d \in \mathbb{R}$ and $0 < \gamma, \eta \leq 1$, then the followings hold:*

1. $\mathscr{L}_x^\eta \mathscr{S}_y^\gamma[c] = \mathscr{L}_x \mathscr{S}_y[c] = \frac{c}{v}, v > 0$
2. $\mathscr{L}_x^\eta \mathscr{S}_y^\gamma\left[\left(\frac{x^\eta}{\eta}\right)^m \left(\frac{y^\gamma}{\gamma}\right)^n\right] = \mathscr{L}_x \mathscr{S}_y[(x)^m (y)^n] = \frac{m!\, n!\, \omega^n}{v^{m+1}}$, where $m$ and $n$ are positive integral.
3. $\mathscr{L}_x^\eta \mathscr{S}_y^\gamma\left[e^{c\frac{x^\eta}{\eta} + d\frac{y^\gamma}{\gamma}}\right] = \mathscr{L}_x \mathscr{S}_y[e^{cx + dy}] = \frac{1}{(v-c)(1-d\omega)}$
4. $\mathscr{L}_x^\eta \mathscr{S}_y^\gamma\left[\sin\left(c\frac{x^\eta}{\eta}\right) \sin\left(d\frac{y^\gamma}{\gamma}\right)\right] = \mathscr{L}_x \mathscr{S}_y[\sin(cx)\sin(dy)] = \frac{c}{(v^2+c^2)} \frac{d\omega}{1+d^2\omega^2}$
5. *First shifting property for CDLST:*

   *If $\mathscr{L}_x^\eta \mathscr{S}_y^\gamma\left[\psi\left(\frac{x^\eta}{\eta}, \frac{y^\gamma}{\gamma}\right)\right] = \Psi_{\eta,\gamma}(v, \omega)$, then*

   $$\mathscr{L}_x^\eta \mathscr{S}_y^\gamma\left[e^{c\frac{x^\eta}{\eta} + d\frac{y^\gamma}{\gamma}} \psi\left(\frac{x^\eta}{\eta}, \frac{y^\gamma}{\gamma}\right)\right] = \frac{1}{1-d\omega} \Psi_{\eta,\gamma}\left(v - c, \frac{\omega}{1-d\omega}\right).$$

6. *Second shifting property for the CDLST:*
   *If $\mathscr{L}_x^\eta \mathscr{S}_y^\gamma\left[\psi\left(\frac{x^\eta}{\eta}, \frac{y^\gamma}{\gamma}\right)\right] = \Psi_{\eta,\gamma}(v, \omega)$, then*

   $$\mathscr{L}_x^\eta \mathscr{S}_y^\gamma\left[\psi\left(\frac{x^\eta}{\eta} - \frac{\delta^\eta}{\eta}, \frac{y^\gamma}{\gamma} - \frac{\varepsilon^\gamma}{\gamma}\right) H\left(\frac{x^\eta}{\eta} - \frac{\delta^\eta}{\eta}, \frac{y^\gamma}{\gamma} - \frac{\varepsilon^\gamma}{\gamma}\right)\right] = e^{-v\frac{\delta^\eta}{\eta} - \frac{1}{\omega}\frac{\varepsilon^\gamma}{\gamma}} \Psi_{\eta,\gamma}(v, \omega),$$

where $H\left(\frac{x^\eta}{\eta} - \frac{\delta^\eta}{\eta}, \frac{y^\gamma}{\gamma} - \frac{\varepsilon^\gamma}{\gamma}\right) = \begin{cases} 1, & \frac{x^\eta}{\eta} > \frac{\delta^\eta}{\eta}, \frac{y^\gamma}{\gamma} > \frac{\varepsilon^\gamma}{\gamma} \\ 0, & \text{otherwise} \end{cases}$

**Proof.** For the proof, see [45–48].

*Now we discuss the conditions for the existence of the CDLST.*

If $\psi$ is of the exponential orders $c$ and; $\left(\frac{x^\eta}{\eta}, \frac{y^\gamma}{\gamma}\right) \to (\infty, \infty)$ and if there exist a nonnegative real number $K : \forall \frac{x^\eta}{\eta} > \mathscr{X}, \frac{y^\gamma}{\gamma} > \mathscr{Y}$; then:

$$\left|\psi\left(\frac{x^\eta}{\eta}, \frac{y^\gamma}{\gamma}\right)\right| = K e^{c\frac{x^\eta}{\eta} + d\frac{y^\gamma}{\gamma}},$$

and we write:

$$\psi\left(\frac{x^\eta}{\eta}, \frac{y^\gamma}{\gamma}\right) = O\left(e^{c\frac{x^\eta}{\eta} + d\frac{y^\gamma}{\gamma}}\right) \text{ as } \left(\frac{x^\eta}{\eta}, \frac{y^\gamma}{\gamma}\right) \to (\infty, \infty)$$

or,

$$\lim_{\left(\frac{x^\eta}{\eta}, \frac{y^\gamma}{\gamma}\right) \to (\infty,\infty)} e^{-v\frac{x^\eta}{\eta} - \frac{1}{\omega}\frac{y^\gamma}{\gamma}} \left|\psi\left(\frac{x^\eta}{\eta}, \frac{y^\gamma}{\gamma}\right)\right| = K \lim_{\left(\frac{x^\eta}{\eta}, \frac{y^\gamma}{\gamma}\right) \to (\infty,\infty)} e^{-(v-c)\frac{x^\eta}{\eta} - \left(\frac{1}{\omega} - d\right)\frac{y^\gamma}{\gamma}} = 0,$$

$$v > c, \frac{1}{\omega} > d.$$

Then, the function $\psi$ is of exponential order as $\left(\frac{x^\eta}{\eta}, \frac{y^\gamma}{\gamma}\right) \to (\infty, \infty)$

**Theorem 8.** Assume that the function $\psi$ defined on the interval $(0, \mathscr{X})$ and $(0, \mathscr{Y})$ of the exponential orders c and d, then CDLST of $\psi$ exists for all $v$ and $\frac{1}{\omega}$ supplied $\text{Re}[v] > c$ and $\text{Re}\left[\frac{1}{\omega}\right] > d$.

**Proof.** We find, from the Definition 6 (iii),

$$\left|\Psi_{\eta,\gamma}(v,\omega)\right| = \left|\frac{1}{\omega}\int_0^\infty\int_0^\infty e^{-v\frac{x^\eta}{\eta} - \frac{1}{\omega}\frac{y^\gamma}{\gamma}}\psi\left(\frac{x^\eta}{\eta}, \frac{y^\gamma}{\gamma}\right)dxdy\right| \le K\int_0^\infty e^{-(v-c)\frac{x^\eta}{\eta}}dx\int_0^\infty \frac{1}{\omega}e^{-\left(\frac{1}{\omega}-d\right)\frac{y^\gamma}{\gamma}}dy = \frac{K}{(v-c)(1-d\omega)}, \text{Re}[v] > c, \text{Re}\left[\frac{1}{\omega}\right] > d \quad (5)$$

Thus, from Equation (5) we've got,

$$\lim_{\left(\frac{x^\eta}{\eta}, \frac{y^\gamma}{\gamma}\right) \to (\infty, \infty)} |\Psi(v,\omega)| = 0, \text{ or } \lim_{\left(\frac{x^\eta}{\eta}, \frac{y^\gamma}{\gamma}\right) \to (\infty, \infty)} \Psi(v,\omega) = 0.$$

**Theorem 9.** If

$$\mathscr{L}_x^\eta \mathscr{S}_y^\gamma\left[\psi\left(\frac{x^\eta}{\eta}, \frac{y^\gamma}{\gamma}\right)\right] = \Psi_{\eta,\gamma}(v,\omega) = \Psi(v,\omega),$$

then the CDLST of the conformable partial derivatives of orders $\eta$ and $\gamma$ can be represented as follows: $0 < \gamma, \eta \le 1$.

$$\mathscr{L}_x^\eta \mathscr{S}_y^\gamma\left[\frac{\partial^\eta}{\partial x^\eta}\left(\psi\left(\frac{x^\eta}{\eta}, \frac{y^\gamma}{\gamma}\right)\right)\right] = v\Psi(v,\omega) - \Psi(0,\omega). \quad (6)$$

$$\mathscr{L}_x^\eta \mathscr{S}_y^\gamma\left[\frac{\partial^\gamma}{\partial y^\gamma}\left(\psi\left(\frac{x^\eta}{\eta}, \frac{y^\gamma}{\gamma}\right)\right)\right] = \omega^{-1}\Psi(v,\omega) - \omega^{-1}\Psi(v,0). \quad (7)$$

$$\mathscr{L}_x^\eta \mathscr{S}_y^\gamma\left[\frac{\partial^{2\eta}}{\partial x^{2\eta}}\left(\psi\left(\frac{x^\eta}{\eta}, \frac{y^\gamma}{\gamma}\right)\right)\right] = v^2\Psi(v,\omega) - v\Psi(0,\omega) - \Psi_x(0,\omega). \quad (8)$$

$$\mathscr{L}_x^\eta \mathscr{S}_y^\gamma\left[\frac{\partial^{2\gamma}}{\partial y^{2\gamma}}\left(\psi\left(\frac{x^\eta}{\eta}, \frac{y^\gamma}{\gamma}\right)\right)\right] = \omega^{-2}\Psi(v,\omega) - \omega^{-2}\Psi(v,0) - \omega^{-1}\Psi_y(v,0). \quad (9)$$

Proof of (6)

$$\mathscr{L}_x^\eta \mathscr{S}_y^\gamma\left[\frac{\partial^\eta}{\partial x^\eta}\left(\psi\left(\frac{x^\eta}{\eta}, \frac{y^\gamma}{\gamma}\right)\right)\right] = \frac{1}{\omega}\int_0^\infty\int_0^\infty e^{-v\frac{x^\eta}{\eta} - \frac{1}{\omega}\frac{y^\gamma}{\gamma}}\frac{\partial^\eta \psi}{\partial x^\eta}x^{\eta-1}y^{\gamma-1}dxdy = \frac{1}{\omega}\int_0^\infty e^{-\frac{1}{\omega}\frac{y^\gamma}{\gamma}}y^{\gamma-1}\left(\int_0^\infty e^{-v\frac{x^\eta}{\eta}}\frac{\partial^\eta \psi}{\partial x^\eta}x^{\eta-1}dx\right)dy. \quad (10)$$

Since we have Theorem 4, $\frac{\partial^\eta \psi}{\partial x^\eta} = x^{-\eta+1}\frac{\partial \psi}{\partial x}$. we use this result into Equation (10).

Therefore, Equation (10) becomes

$$\mathscr{L}_x^\eta \mathscr{S}_y^\gamma\left[\frac{\partial^\eta}{\partial x^\eta}\left(\psi\left(\frac{x^\eta}{\eta}, \frac{y^\gamma}{\gamma}\right)\right)\right] = \frac{1}{\omega}\int_0^\infty\int_0^\infty e^{-v\frac{x^\eta}{\eta} - \frac{1}{\omega}\frac{y^\gamma}{\gamma}}\frac{\partial^\eta \psi}{\partial x^\eta}x^{\eta-1}y^{\gamma-1}dxdy = \frac{1}{\omega}\int_0^\infty e^{-\frac{1}{\omega}\frac{y^\gamma}{\gamma}}y^{\gamma-1}\left(\int_0^\infty e^{-v\frac{x^\eta}{\eta}}\frac{\partial \psi}{\partial x}dx\right)dy. \quad (11)$$

*Considering*

$$\int_0^\infty e^{-v\frac{x^\eta}{\eta}}\frac{\partial \psi}{\partial x}dx = -\psi(0,y) + v\Psi(v,y) \quad (12)$$

*Substituting* (12) *in* (11) *with simple calculations, we get the result* (6), *as follows*

$$\mathscr{L}_x^\eta \mathscr{S}_y^\gamma\left[\frac{\partial^\eta}{\partial y^\eta}\left(\psi\left(\frac{x^\eta}{\eta}, \frac{y^\gamma}{\gamma}\right)\right)\right] = v\Psi(v,\omega) - \Psi(0,\omega).$$

*The remaining results in* (7), (8), *and* (9) *can be proved by similar arguments.*
*Additionally, the previous findings can be generalized I the following formulas.*

$$\mathscr{L}_x^\eta \mathscr{S}_y^\gamma \left[ \frac{\partial^{m\eta}}{\partial x^{m\eta}} \left( \psi \left( \frac{x^\eta}{\eta}, \frac{y^\gamma}{\gamma} \right) \right) \right] = v^m \Psi(v, \omega) - \sum_{k=0}^{m-1} v^{m-1-k} \mathscr{S}_y^\gamma \left[ \frac{\partial^{k\eta}}{\partial x^{k\eta}} \left( \psi \left( 0, \frac{y^\gamma}{\gamma} \right) \right) \right], \quad (13)$$

$$\mathscr{L}_x^\eta \mathscr{S}_y^\gamma \left[ \frac{\partial^{n\gamma}}{\partial y^{n\gamma}} \left( \psi \left( \frac{x^\eta}{\eta}, \frac{y^\gamma}{\gamma} \right) \right) \right] = \omega^{-n} \Psi(v, \omega) - \sum_{j=0}^{n-1} \omega^{-n+j} \mathscr{L}_x^\eta \left[ \frac{\partial^{j\gamma}}{\partial y^{j\gamma}} \left( \psi \left( \frac{x^\eta}{\eta}, 0 \right) \right) \right]. \quad (14)$$

*The proof of* (13) *and* (14) *can be obtained by mathematical induction.*

## 4. Idea of CDLST combined with modified decomposition technique

In this part, we exercised CDLST joined to the modified decomposition approach to find the solutions to nonlinear conformable fractional equations as below.

Applying CDLST on Equation (1), we get

$$b \left( \omega^{-n} \Psi(v, \omega) - \sum_{j=0}^{n-1} \omega^{-n+j} \mathscr{L}_x^\eta \left[ \frac{\partial^{j\gamma} \psi(\frac{x^\eta}{\eta}, 0)}{\partial y^{j\gamma}} \right] \right) + c \left( v^m \Psi(v, \omega) - \sum_{k=0}^{m-1} v^{m-1-k} \mathscr{S}_y^\gamma \left[ \frac{\partial^{k\eta} \psi(0, \frac{y^\gamma}{\gamma})}{\partial y^{k\eta}} \right] \right) + \mathscr{L}_x^\eta \mathscr{S}_y^\gamma \left[ \mathcal{N} \left[ \psi \left( \frac{x^\eta}{\eta}, \frac{y^\gamma}{\gamma} \right) \right] \right] = G_1(v, \omega)$$
$$+ \mathscr{L}_x^\eta \mathscr{S}_y^\gamma \left[ g_2 \left( \frac{x^\eta}{\eta}, \frac{y^\gamma}{\gamma} \right) \right]. \quad (15)$$

Operating the CLT and CST for conditions (2) and (3), we get,

$$\mathscr{L}_x^\eta \left[ \frac{\partial^{j\gamma} \psi(\frac{x^\eta}{\eta}, 0)}{\partial y^{j\gamma}} \right] = F_j(v), j = 0, 1, \ldots, n-1,$$

$$\mathscr{S}_y^\gamma \left[ \frac{\partial^{k\eta} \psi(0, \frac{y^\gamma}{\gamma})}{\partial x^{k\eta}} \right] = H_k(\omega), k = 0, 1, \ldots, m-1. \quad (16)$$

By substituting (16) for (15), we have

$$b \left( \omega^{-n} \Psi(v, \omega) - \sum_{j=0}^{n-1} \omega^{-n+j} F_j(v) \right) + c \left( v^m \Psi(v, \omega) - \sum_{k=0}^{m-1} v^{m-1-k} H_k(\omega) \right) = G_1(v, \omega) + \mathscr{L}_x^\eta \mathscr{S}_y^\gamma \left[ g_2 \left( \frac{x^\eta}{\eta}, \frac{y^\gamma}{\gamma} \right) - \mathcal{N} \left[ \psi \left( \frac{x^\eta}{\eta}, \frac{y^\gamma}{\gamma} \right) \right] \right]. \quad (17)$$

Simplifying Equation (17), we obtain

$$\Psi(v, \omega) = [b\omega^{-n} + cv^m]^{-1} \left( b \left( \sum_{j=0}^{n-1} \omega^{-n+j} F_j(v) \right) + c \left( \sum_{k=0}^{m-1} v^{m-1-k} H_k(\omega) \right) + G_1(v, \omega) \right)$$
$$+ (b\omega^{-n} + cv^m)^{-1} \mathscr{L}_x^\eta \mathscr{S}_y^\gamma \left[ g_2 \left( \frac{x^\eta}{\eta}, \frac{y^\gamma}{\gamma} \right) - \mathcal{N} \left[ \psi \left( \frac{x^\eta}{\eta}, \frac{y^\gamma}{\gamma} \right) \right] \right] \quad (18)$$

Taking $(\mathscr{L}_x^\eta)^{-1}(\mathscr{S}_y^\gamma)^{-1}$ of (18), we get

$$\psi \left( \frac{x^\eta}{\eta}, \frac{y^\gamma}{\gamma} \right) = (\mathscr{L}_x^\eta)^{-1} (\mathscr{S}_y^\gamma)^{-1} \left[ [b\omega^{-n} + cv^m]^{-1} \left( b \left( \sum_{j=0}^{n-1} \omega^{-n+j} F_j(v) \right) + c \left( \sum_{k=0}^{m-1} v^{m-1-k} H_k(\omega) \right) + G_1(v, \omega) \right) \right]$$
$$+ (\mathscr{L}_x^\eta)^{-1} (\mathscr{S}_y^\gamma)^{-1} \left[ (b\omega^{-n} + cv^m)^{-1} \mathscr{L}_x^\eta \mathscr{S}_y^\gamma \left[ g_2 \left( \frac{x^\eta}{\eta}, \frac{y^\gamma}{\gamma} \right) - \mathcal{N} \left[ \psi \left( \frac{x^\eta}{\eta}, \frac{y^\gamma}{\gamma} \right) \right] \right] \right]. \quad (19)$$

Now, use the modified decomposition approach by assuming

$$\psi \left( \frac{x^\eta}{\eta}, \frac{y^\gamma}{\gamma} \right) = \sum_{i=0}^{\infty} \psi_i \left( \frac{x^\eta}{\eta}, \frac{y^\gamma}{\gamma} \right) \quad (20)$$

By decomposing the nonlinear term as:

$$\mathcal{N} \left[ \psi \left( \frac{x^\eta}{\eta}, \frac{y^\gamma}{\gamma} \right) \right] = \sum_{i=0}^{\infty} A_i, \quad (21)$$

for some Adomian polynomials $A_i(\psi)$ that are given by:

$$A_i(\psi_0, \psi_1, \psi_2, ..., \psi_n) = \frac{1}{i!} \frac{d^i}{d\lambda^i}\left(\mathcal{N}\left[\sum_{i=0}^{\infty} \lambda^i \psi_i\right]\right)_{\lambda=0}, i = 0, 1, 2, \cdots.$$

Substituting Equation (20) and Equation (21) in Equation (19), we get:

$$\sum_{i=0}^{\infty} \psi_i\left(\frac{x^\eta}{\eta}, \frac{y^\gamma}{\gamma}\right) = (\mathcal{L}_x^\eta)^{-1}(\mathcal{S}_y^\gamma)^{-1}\left[(b\omega^{-n} + cv^m)^{-1}\left(b\left(\sum_{j=0}^{n-1} \omega^{-n+j} F_j(v)\right) + c\left(\sum_{k=0}^{m-1} v^{m-1-k} H_k(\omega)\right) + G_1(v, \omega)\right)\right.$$

$$\left. + (\mathcal{L}_x^\eta)^{-1}(\mathcal{S}_y^\gamma)^{-1}\left[[b\omega^{-n} + cv^m]^{-1}\mathcal{L}_x^\eta \mathcal{S}_y^\gamma\left[g_2\left(\frac{x^\eta}{\eta}, \frac{y^\gamma}{\gamma}\right) - \sum_{i=0}^{\infty} A_i\right]\right]\right].$$

As result, the following recurrence relations are obtained

$$\psi_0\left(\frac{x^\eta}{\eta}, \frac{y^\gamma}{\gamma}\right) = (\mathcal{L}_x^\eta)^{-1}(\mathcal{S}_y^\gamma)^{-1}\left[(b\omega^{-n} + cv^m)^{-1}\left(b\left(\sum_{j=0}^{n-1} \omega^{-n+j} F_j(v)\right) + c\left(\sum_{k=0}^{m-1} v^{m-1-k} H_k(\omega)\right) + G_1(v, \omega)\right)\right], \quad (22)$$

$$\psi_1\left(\frac{x^\eta}{\eta}, \frac{y^\gamma}{\gamma}\right) = (\mathcal{L}_x^\eta)^{-1}(\mathcal{S}_y^\gamma)^{-1}\left[(b\omega^{-n} + cv^m)^{-1} \mathcal{L}_x^\eta \mathcal{S}_y^\gamma\left[g_2\left(\frac{x^\eta}{\eta}, \frac{y^\gamma}{\gamma}\right) - A_0\right]\right], \quad (23)$$

$$\psi_{s+1}\left(\frac{x^\eta}{\eta}, \frac{y^\gamma}{\gamma}\right) = -(\mathcal{L}_x^\eta)^{-1}(\mathcal{S}_y^\gamma)^{-1}\left[(b\omega^{-n} + cv^m)^{-1} \mathcal{L}_x^\eta \mathcal{S}_y^\gamma[A_s]\right], s \geq 1. \quad (24)$$

Following that, one can get the solution of (1)

$$\psi\left(\frac{x^\eta}{\eta}, \frac{y^\gamma}{\gamma}\right) = \psi_0\left(\frac{x^\eta}{\eta}, \frac{y^\gamma}{\gamma}\right) + \psi_1\left(\frac{x^\eta}{\eta}, \frac{y^\gamma}{\gamma}\right) + \psi_2\left(\frac{x^\eta}{\eta}, \frac{y^\gamma}{\gamma}\right) + \cdots.$$

## 5. Numerical applications

To assess the effectiveness and usefulness of the CDLSMD approach, numerous applications of nonlinear fractional partial differential equations are discussed in this section.

**Example 4.1.** Let us dicuss the dissipative wave equation of conformable derivative as:

$$\frac{\partial^{2\gamma} \psi}{\partial y^{2\gamma}} - \frac{\partial^2 \psi}{\partial x^2} + \frac{\partial}{\partial y}\left(\psi \frac{\partial \psi}{\partial x}\right) = 2e^{-\frac{y^\gamma}{\gamma}} \sin\frac{x^\eta}{\eta} - 2e^{-2\frac{y^\gamma}{\gamma}} \sin\frac{x^\eta}{\eta}\cos\frac{x^\eta}{\eta}, 0 < \gamma \leq 1 \quad (25)$$

with initial conditions:

$$\psi\left(\frac{x^\eta}{\eta}, 0\right) = \sin\frac{x^\eta}{\eta}, \psi_y\left(\frac{x^\eta}{\eta}, 0\right) = -\sin\frac{x^\eta}{\eta}, \quad (26)$$

and BCs:

$$\psi\left(0, \frac{y^\gamma}{\gamma}\right) = 0, \psi_x\left(0, \frac{y^\gamma}{\gamma}\right) = e^{-\frac{y^\gamma}{\gamma}}. \quad (27)$$

**Solution.** Using the CDLST to (25) and the CLT to (26) and CST to (27), to get,

$$\Psi(v, \omega) = \frac{1}{(v^2 + 1)(1 + \omega)} - \frac{\omega^2}{1 - v^2\omega^2}\mathcal{L}_x^\eta \mathcal{S}_y^\gamma\left[2e^{-2\frac{y^\gamma}{\gamma}} \sin\frac{x^\eta}{\eta}\cos\frac{x^\eta}{\eta} + \frac{\partial}{\partial y}\left(\psi \frac{\partial \psi}{\partial x}\right)\right]. \quad (28)$$

Taking $(\mathcal{L}_x^\eta)^{-1}(\mathcal{S}_y^\gamma)^{-1}[\Psi(v, \omega)]$ of (28), we get

$$\psi\left(\frac{x^\eta}{\eta}, \frac{y^\gamma}{\gamma}\right) = e^{-\frac{y^\gamma}{\gamma}} \sin\frac{x^\eta}{\eta} - (\mathcal{L}_x^\eta)^{-1}(\mathcal{S}_y^\gamma)^{-1}\left[\frac{\omega^2}{1 - v^2\omega^2}\mathcal{L}_x^\eta \mathcal{S}_y^\gamma\left[2e^{-2\frac{y^\gamma}{\gamma}} \sin\frac{x^\eta}{\eta}\cos\frac{x^\eta}{\eta} + \frac{\partial}{\partial y}\left(\psi \frac{\partial \psi}{\partial x}\right)\right]\right]. \quad (29)$$

Following that, operating the modified decomposition approach and substituting (20) in (29) and utilizing the outcomes in (22), (23) and (24), we get the components of the solution as:

$$\psi_0\left(\frac{x^\eta}{\eta}, \frac{y^\gamma}{\gamma}\right) = e^{-\frac{y^\gamma}{\gamma}} \sin\frac{x^\eta}{\eta},$$

$$\psi_1\left(\frac{x^\eta}{\eta}, \frac{y^\gamma}{\gamma}\right) = -(\mathcal{L}_x^\eta)^{-1}(\mathcal{S}_y^\gamma)^{-1}\left[\frac{\omega^2}{1 - v^2\omega^2}\mathcal{L}_x^\eta \mathcal{S}_y^\gamma\left[2e^{-2\frac{y^\gamma}{\gamma}} \sin\frac{x^\eta}{\eta}\cos\frac{x^\eta}{\eta} + A_0\right]\right] = 0, A_0 = \frac{\partial}{\partial y}(\psi_0(\psi_0)_x),$$

$$\psi_2\left(\frac{x^\eta}{\eta},\frac{y^\gamma}{\gamma}\right) = -\left(\mathscr{L}_x^\eta\right)^{-1}\left(\mathscr{S}_y^\gamma\right)^{-1}\left[\frac{\omega^2}{1-v^2\omega^2}\mathscr{L}_x^\eta\mathscr{S}_y^\gamma[A_1]\right] = 0, A_1 = \frac{\partial}{\partial y}(\psi_1\psi_{0_x}+\psi_0\psi_{1_x}).$$

Thus, we obtain the accurate solution of Equation (25) as:

$$\psi\left(\frac{x^\eta}{\eta},\frac{y^\gamma}{\gamma}\right) = e^{-\frac{y^\gamma}{\gamma}}\sin\frac{x^\eta}{\eta}, \qquad (30)$$

and putting $\eta = \gamma = 1$; we have the exact solution as,

$$\psi(x,y) = e^{-y}\sin x.$$

In the following figures, we sketch the 3D plots of the accurate solution of Equation (25) in Fig. 1 (a), which is the same approximate obtained solution, when putting $\eta = \gamma = 1$. In Fig. 1 (b), we sketch the approximate solution with various values of $\eta$ and $\gamma$: $\eta = \gamma = 1, 0.9, 0.8$.

In Table 1 below, we introduce the absolute exact errors obtained from computing the absolute difference of the exact and CDLSMD solutions obtained in this example with different values of the variables $y = 0.1, 0.2, 0.3$ at $x = 1$.

**Example 4.2.** The KdV equation of conformable derivative is given by:

$$\frac{\partial^\gamma\psi}{\partial y^\gamma} - \psi\frac{\partial\psi}{\partial x} + \frac{\partial^3\psi}{\partial x^3} = -e^{\frac{x^\eta}{\eta}}\left(\frac{y^\gamma}{\gamma}+1\right) + \frac{y^\gamma}{\gamma}e^{\frac{x^\eta}{\eta}}\left(-\frac{y^\gamma}{\gamma}e^{\frac{x^\eta}{\eta}}+1\right), 0 < \gamma \le 1 \qquad (31)$$

with initial condition:

$$\psi\left(\frac{x^\eta}{\eta},0\right) = 1. \qquad (32)$$

and boundary conditions:

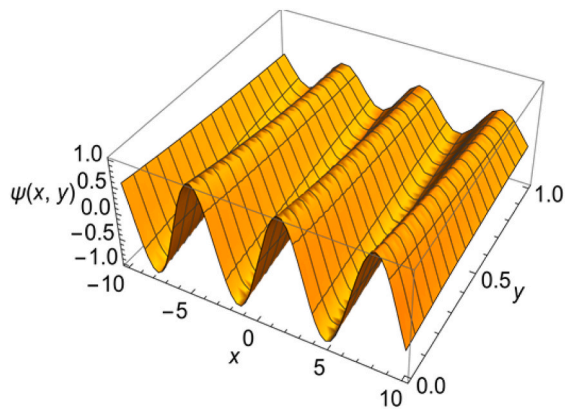

(a)

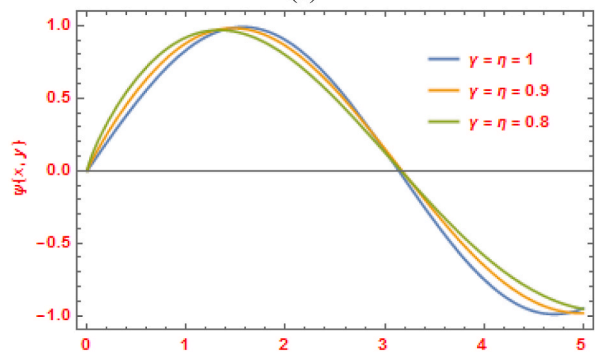

(b)

**Fig. 1.** (a) The 3D plots of the solution of Equation (30) gained by the presented method comparing to the exact solution, (b) The CDLSMD solution of $\psi\left(\frac{x^\eta}{\eta},\frac{y^\gamma}{\gamma}\right)$ for Equation (30) at $\eta = \gamma = 1, 0.9, 0.8$.

**Table 1**

The absolute error of $\psi\left(\frac{x^\eta}{\eta},\frac{y^\gamma}{\gamma}\right)$ given by the CDLSMD technique for example 4.1 at various values of $\eta$, $\gamma$ and $y$.

| x | y | η | γ | Exact | CDLSMD | $\|\psi_{Exact} - \psi_{CDLSMD}\|$ |
|---|---|---|---|---|---|---|
| 1 | 0.1 | 0.8 | 0.8 | 0.761394 | 0.778431 | $1.7037 \times 10^{-2}$ |
|   | 0.2 |     |     | 0.688938 | 0.672136 | $1.68025 \times 10^{-2}$ |
|   | 0.3 |     |     | 0.623377 | 0.588923 | $3.44542 \times 10^{-2}$ |
| 1 | 0.1 | 0.9 | 0.9 | 0.761394 | 0.779205 | $1.78107 \times 10^{-2}$ |
|   | 0.2 |     |     | 0.688938 | 0.690302 | $1.36333 \times 10^{-3}$ |
|   | 0.3 |     |     | 0.623377 | 0.615339 | $8.03807 \times 10^{-3}$ |

$$\psi\left(0,\frac{y^\gamma}{\gamma}\right) = 1 - \frac{y^\gamma}{\gamma}, \psi_x\left(0,\frac{y^\gamma}{\gamma}\right) = \psi_{xx}\left(0,\frac{y^\gamma}{\gamma}\right) = -\frac{y^\gamma}{\gamma}. \tag{33}$$

**Solution.** Applying the CDLST on (31) and CLT to (32) and the CST to (33), to get

$$\Psi(v,\omega) = \frac{1}{v} - \frac{\omega}{v-1} + \frac{\omega}{1+\omega v^3}\mathscr{L}_x^\eta \mathscr{S}_y^\gamma\left[\frac{y^\gamma}{\gamma}e^{\frac{x^\eta}{\eta}}\left(-\frac{y^\gamma}{\gamma}e^{\frac{x^\eta}{\eta}} + 1\right) + \psi\frac{\partial\psi}{\partial x}\right]. \tag{34}$$

Taking $(\mathscr{L}_x^\eta)^{-1}(\mathscr{S}_y^\gamma)^{-1}[\Psi(v,\omega)]$ of (34), we get

$$\psi\left(\frac{x^\eta}{\eta},\frac{y^\gamma}{\gamma}\right) = 1 - \frac{y^\gamma}{\gamma}e^{\frac{x^\eta}{\eta}} + (\mathscr{L}_x^\eta)^{-1}(\mathscr{S}_y^\gamma)^{-1}\left[\frac{\omega}{1+\omega v^3}\mathscr{L}_x^\eta \mathscr{S}_y^\gamma\left[\frac{y^\gamma}{\gamma}e^{\frac{x^\eta}{\eta}}\left(-\frac{y^\gamma}{\gamma}e^{\frac{x^\eta}{\eta}} + 1\right) + \psi\frac{\partial\psi}{\partial x}\right]\right]. \tag{35}$$

Using the modified decomposition technique, and substituting (20) in (35) and with the findings in Equations (22)–(24), we get the solution components

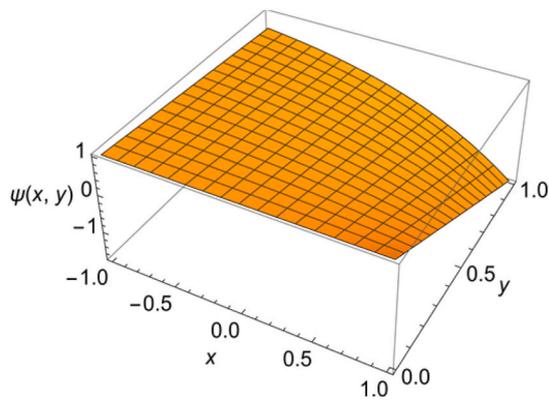

(a)

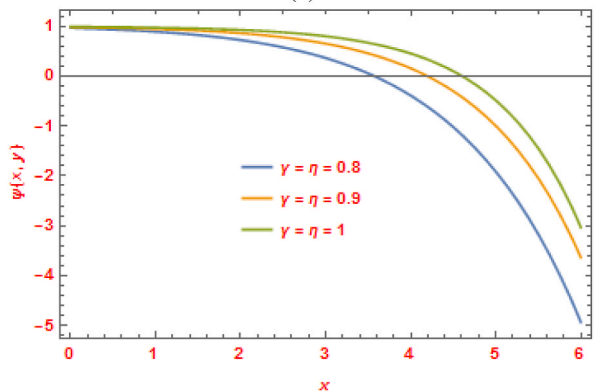

(b)

**Fig. 2.** (a) The 3D plots solution of Equation (36) gained using the presented method comparing to the exact solution (b) The CDLSMD solution of $\psi\left(\frac{x^\eta}{\eta},\frac{y^\gamma}{\gamma}\right)$ for Equation (36) at $\eta = \gamma = 1, 0.9, 0.8$.

$$\psi_0\left(\frac{x^\eta}{\eta},\frac{y^\gamma}{\gamma}\right) = 1 - \frac{y^\gamma}{\gamma}e^{\frac{x^\eta}{\eta}},$$

$$\psi_1\left(\frac{x^\eta}{\eta},\frac{y^\gamma}{\gamma}\right) = (\mathscr{L}_x^\eta)^{-1}(\mathscr{S}_y^\gamma)^{-1}\left[\frac{\omega}{1+\omega v^3}\mathscr{L}_x^\eta\mathscr{S}_y^\gamma\left[\frac{y^\gamma}{\gamma}e^{\frac{x^\eta}{\eta}}\left(-\frac{y^\gamma}{\gamma}e^{\frac{x^\eta}{\eta}}+1\right)+A_0\right]\right] = 0, A_0 = \psi_0\frac{\partial\psi_0}{\partial x},$$

$$\psi_2\left(\frac{x^\eta}{\eta},\frac{y^\gamma}{\gamma}\right) = (\mathscr{L}_x^\eta)^{-1}(\mathscr{S}_y^\gamma)^{-1}\left[\frac{\omega}{1+\omega v^3}\mathscr{L}_x^\eta\mathscr{S}_y^\gamma[A_1]\right] = 0, A_1 = \psi_1\psi_{0_x} + \psi_0\psi_{1_x}.$$

Thus, we get the solution to the (31) in the form:

$$\psi\left(\frac{x^\eta}{\eta},\frac{y^\gamma}{\gamma}\right) = 1 - \frac{y^\gamma}{\gamma}e^{\frac{x^\eta}{\eta}}. \tag{36}$$

Note that putting $\eta = \gamma = 1$; then the exact solution is,

$$\psi(x,y) = 1 - ye^x.$$

In the following figures, we sketch the 3D plot of the accurate solution of Example 4.2 in Fig. 2(a) which is the same approximate obtained solution, when putting $\eta = \gamma = 1$. In Fig. 2(b), we sketch the approximate solution with different values of the fractional orders $\eta$ and $\gamma$: $\eta = \gamma = 1, 0.9, 0.8$.

In Table 2 below, we introduce the absolute exact errors obtained from computing the absolute difference of the accurate and CDLSMD solutions obtained in this application with different values of the variables $y = 0.1, 0.2, 0.3$ at $x = 1$.

**Example 4.3.** The nonlinear advection equation of conformable partial differential equation:

$$\frac{\partial^\gamma\psi}{\partial y^\gamma} + \psi\frac{\partial\psi}{\partial x} = 2\frac{y^\gamma}{\gamma} + \frac{x^\eta}{\eta} + \left(\frac{y^\gamma}{\gamma}\right)^3 + \frac{x^\eta}{\eta}\left(\frac{y^\gamma}{\gamma}\right)^2, 0 < \gamma \le 1, \tag{37}$$

with initial condition:

$$\psi\left(\frac{x^\eta}{\eta},0\right) = 0. \tag{38}$$

**Solution.** Operating the CDLST to Equation (37) and the CLT to Equation (38), to get

$$\Psi(v,\omega) = \frac{2\omega^2}{v} + \frac{\omega}{v^2} + \omega\mathscr{L}_x^\eta\mathscr{S}_y^\gamma\left[\left(\frac{y^\gamma}{\gamma}\right)^3 + \frac{x^\eta}{\eta}\left(\frac{y^\gamma}{\gamma}\right)^2 - \psi\frac{\partial\psi}{\partial x}\right]. \tag{39}$$

Taking $(\mathscr{L}_x^\eta)^{-1}(\mathscr{S}_y^\gamma)^{-1}[\Psi(v,\omega)]$ of (39), we get

$$\psi\left(\frac{x^\eta}{\eta},\frac{y^\gamma}{\gamma}\right) = \left(\frac{y^\gamma}{\gamma}\right)^2 + \frac{x^\eta}{\eta}\frac{y^\gamma}{\gamma} + (\mathscr{L}_x^\eta)^{-1}(\mathscr{S}_y^\gamma)^{-1}\left[\omega\mathscr{L}_x^\eta\mathscr{S}_y^\gamma\left[\left(\frac{y^\gamma}{\gamma}\right)^3 + \frac{x^\eta}{\eta}\left(\frac{y^\gamma}{\gamma}\right)^2 - \psi\frac{\partial\psi}{\partial x}\right]\right]. \tag{40}$$

Now, applying the modified decomposition method and substitute (20) in (40) and using results in (22), (23), and (24), we get the following solution components

$$\psi_0\left(\frac{x^\eta}{\eta},\frac{y^\gamma}{\gamma}\right) = \left(\frac{y^\gamma}{\gamma}\right)^2 + \frac{x^\eta}{\eta}\frac{y^\gamma}{\gamma},$$

$$\psi_1\left(\frac{x^\eta}{\eta},\frac{y^\gamma}{\gamma}\right) = (\mathscr{L}_x^\eta)^{-1}(\mathscr{S}_y^\gamma)^{-1}\left[\omega\mathscr{L}_x^\eta\mathscr{S}_y^\gamma\left[\left(\frac{y^\gamma}{\gamma}\right)^3 + \frac{x^\eta}{\eta}\left(\frac{y^\gamma}{\gamma}\right)^2 - A_0\right]\right] = 0, A_0 = \psi_0\frac{\partial\psi_0}{\partial x},$$

**Table 2**

The absolute error of $\psi\left(\frac{x^\eta}{\eta},\frac{y^\gamma}{\gamma}\right)$ gained by the CDLSMD technique for Equation (31) at different values of $\eta$, $\gamma$ and $y$.

| x | y | $\eta$ | $\gamma$ | Exact | CDLSMD | $|\psi_{Exact} - \psi_{CDLSMD}|$ |
|---|---|---|---|---|---|---|
| 1 | 0.1 | 0.8 | 0.8 | 0.728172 | 0.308522 | $4.19649 \times 10^{-1}$ |
|   | 0.2 |     |     | 0.456344 | − 0.203932 | $2.52411 \times 10^{-1}$ |
|   | 0.3 |     |     | 0.184515 | − 0.665233 | $8.49749 \times 10^{-1}$ |
| 1 | 0.1 | 0.9 | 0.9 | 0.728172 | 0.57508 | $1.53092 \times 10^{-1}$ |
|   | 0.2 |     |     | 0.456344 | 0.207072 | $2.49272 \times 10^{-1}$ |
|   | 0.3 |     |     | 0.184515 | − 0.14213 | $3.26647 \times 10^{-1}$ |

$$\psi_2\left(\frac{x^\eta}{\eta},\frac{y^\gamma}{\gamma}\right) = \left(\mathcal{L}^\eta_x\right)^{-1}\left(\mathcal{S}^\gamma_y\right)^{-1}\left[\omega\,\mathcal{L}^\eta_x\mathcal{S}^\gamma_y[A_1]\right] = 0,\ A_1 = \psi_1\psi_{0_y} + \psi_0\psi_{1_y}.$$

Following that, we get the solution of Equation (37):

$$\psi\left(\frac{x^\eta}{\eta},\frac{y^\gamma}{\gamma}\right) = \left(\frac{y^\gamma}{\gamma}\right)^2 + \frac{x^\eta}{\eta}\frac{y^\gamma}{\gamma}. \qquad (41)$$

Putting $\eta = \gamma = 1$; one can obtain the solution,

$$\psi(x,y) = y^2 + yx.$$

In the following figures, we sketch the 3D plot of the accurate solution of Equation (37) in Fig. 3 (a) which is the same approximate obtained solution, when putting $\eta = \gamma = 1$. Fig. 3 (b), present the graph of the obtained solution with different values of the fractional orders $\eta$ and $\gamma$: $\eta = \gamma = 1, 0.9, 0.8$.

In Table 3 below, we present the absolute exact errors obtained from computing the absolute difference of the exact and CDLSMD solutions obtained in this problem with different values of the variables $y = 0.1, 0.2, 0.3$ at $x = 1$.

**Example 4.4.** The nonlinear Klein-Gordon problem of conformable partial differential equation:

$$\frac{\partial^{2\gamma}\psi}{\partial y^{2\gamma}} - \frac{\partial^2\psi}{\partial x^2} + \psi^2 = 2\left(\frac{x^\eta}{\eta}\right)^2 - 2\left(\frac{y^\gamma}{\gamma}\right)^2 + \left(\frac{x^\eta}{\eta}\right)^4\left(\frac{y^\gamma}{\gamma}\right)^4,\ 0 < \gamma \leq 1, \qquad (42)$$

with initial conditions:

$$\psi\left(\frac{x^\eta}{\eta},0\right) = 0,\ \psi_y\left(\frac{x^\eta}{\eta},0\right) = 0, \qquad (43)$$

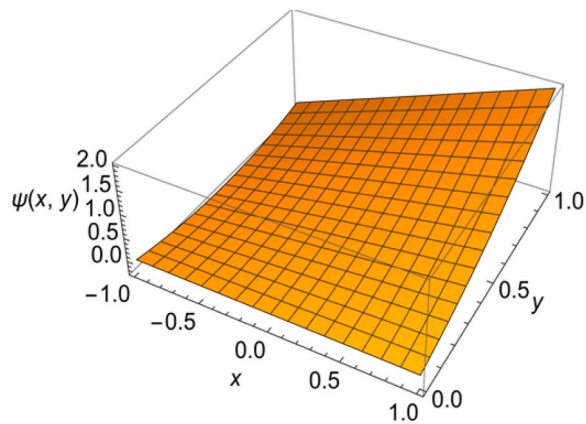

(a)

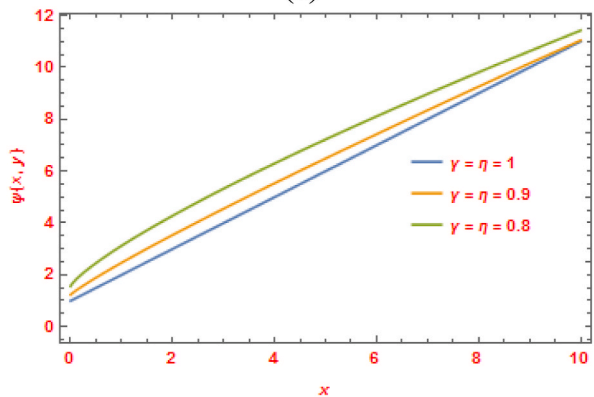

(b)

**Fig. 3.** (a) The 3D plots of the solution graph of Equation (41) gained by the presented technique comparing with accurate solution (b) The CLSMD solution of $\psi\left(\frac{x^\eta}{\eta},\frac{y^\gamma}{\gamma}\right)$ for Equation (41) at $\eta = \gamma = 1, 0.9, 0.8$.

**Table 3**
The absolute error of $\psi\left(\frac{x^\eta}{\eta},\frac{y^\gamma}{\gamma}\right)$ obtained by CDLSMD method for Equation (37) at various values of $\eta$, $\gamma$ and $y$.

| x | y | η | γ | Exact | **CDLSMD** | $\|\psi_{Exact} - \psi_{CDLSMD}\|$ |
|---|---|---|---|---|---|---|
| 1 | 0.1 | 0.8 | 0.8 | 0.11 | 0.286888 | $1.76888 \times 10^{-1}$ |
|   | 0.2 |     |     | 0.24 | 0.55014  | $3.10144 \times 10^{-1}$ |
|   | 0.3 |     |     | 0.39 | 0.823994 | $4.33994 \times 10^{-1}$ |
| 1 | 0.1 | 0.9 | 0.9 | 0.11 | 0.174989 | $2.15011 \times 10^{-1}$ |
|   | 0.2 |     |     | 0.24 | 0.358164 | $1.18164 \times 10^{-1}$ |
|   | 0.3 |     |     | 0.39 | 0.55912  | $1.6912 \times 10^{-1}$ |

and boundary conditions:

$$\psi\left(0,\frac{y^\gamma}{\gamma}\right) = 0, \psi_x\left(0,\frac{y^\gamma}{\gamma}\right) = 0. \quad (44)$$

**Solution.** Applying the CDLST to (42) and the CLT to (43) and the CST to (44), one can get

$$\Psi(v,\omega) = \frac{4\omega^2}{v^3} + \frac{\omega^2}{1 - v^2\omega^2}\mathscr{L}_x^\eta\mathscr{S}_y^\gamma\left[\left(\frac{x^\eta}{\eta}\right)^4\left(\frac{y^\gamma}{\gamma}\right)^4 - \psi^2\right]. \quad (45)$$

Taking $(\mathscr{L}_x^\eta)^{-1}(\mathscr{S}_y^\gamma)^{-1}[\Psi(v,\omega)]$ of (45), we get

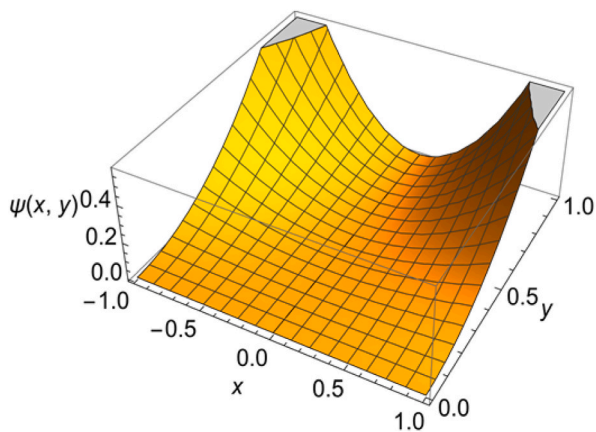

(a)

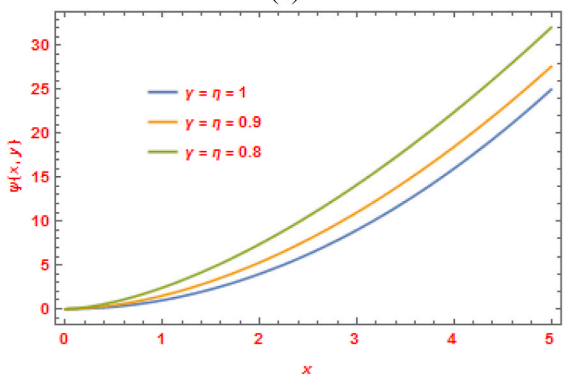

(b)

**Fig. 4.** (a) The 3D solution graphs of Equation (47) obtained by the presented method comparison with exact solution (b) The CDLSMD solution of $\psi\left(\frac{x^\eta}{\eta},\frac{y^\gamma}{\gamma}\right)$ for Equation (47) at $\eta = \gamma = 1, 0.9, 0.8$.

$$\psi\left(\frac{x^\eta}{\eta},\frac{y^\gamma}{\gamma}\right)=\left(\frac{x^\eta}{\eta}\right)^2\left(\frac{y^\gamma}{\gamma}\right)^2+\left(\mathscr{L}_x^\eta\right)^{-1}\left(\mathscr{S}_y^\gamma\right)^{-1}\left[\frac{\omega^2}{1-v^2\omega^2}\mathscr{L}_x^\eta\mathscr{S}_y^\gamma\left[\left(\frac{x^\eta}{\eta}\right)^4\left(\frac{y^\gamma}{\gamma}\right)^4-\psi^2\right]\right]. \tag{46}$$

Using the modified decomposition method, and substituting (20) in (46) and using the outcomes in (22), (23), and (24), we get the components of the solution as

$$\psi_0\left(\frac{x^\eta}{\eta},\frac{y^\gamma}{\gamma}\right)=\left(\frac{x^\eta}{\eta}\right)^2\left(\frac{y^\gamma}{\gamma}\right)^2,$$

$$\psi_1\left(\frac{x^\eta}{\eta},\frac{y^\gamma}{\gamma}\right)=\left(\mathscr{L}_x^\eta\right)^{-1}\left(\mathscr{S}_y^\gamma\right)^{-1}\left[\frac{\omega^2}{1-v^2\omega^2}\mathscr{L}_x^\eta\mathscr{S}_y^\gamma\left[\left(\frac{x^\eta}{\eta}\right)^4\left(\frac{y^\gamma}{\gamma}\right)^4-A_0\right]\right]=0, A_0=\psi_0^2,$$

$$\psi_2\left(\frac{x^\eta}{\eta},\frac{y^\gamma}{\gamma}\right)=\left(\mathscr{L}_x^\eta\right)^{-1}\left(\mathscr{S}_y^\gamma\right)^{-1}\left[\frac{\omega^2}{1-v^2\omega^2}\mathscr{L}_x^\eta\mathscr{S}_y^\gamma[A_1]\right]=0, A_1=2\psi_0\psi_1.$$

Thus, we get the solution of (42) in the form

$$\psi\left(\frac{x^\eta}{\eta},\frac{y^\gamma}{\gamma}\right)=\left(\frac{x^\eta}{\eta}\right)^2\left(\frac{y^\gamma}{\gamma}\right)^2, \tag{47}$$

Putting $\eta=\gamma=1$; then the solution is,

$$\psi(x,x)=x^2y^2.$$

In the figures below, we sketch the 3D plot of the accurate solution of Equation (42) in Fig. 4 (a) which is the same approximate obtained solution, when putting $\eta=\gamma=1$. In Fig. 4 (b), we sketch the approximate solution at various values of the fractional orders $\eta$ and $\gamma$: $\eta=\gamma=1,0.9,0.8$.

In the following table (Table 4), we introduce the absolute exact errors obtained from computing the absolute difference of the exact and CDLSMD solutions obtained in this example with different values of the variables $y=0.1,0.2,0.3$ at $x=1$.

**Example 4.5.** The telegraph equation of nonlinear conformable partial differential equation is given by:

$$\frac{\partial^{2\eta}\psi}{\partial x^{2\eta}}=\frac{\partial^{2\gamma}\psi}{\partial y^{2\gamma}}+2\frac{\partial^\gamma\psi}{\partial y^\gamma}+\psi^2+e^{\frac{x^\eta}{\eta}-2\frac{y^\gamma}{\gamma}}-e^{2\frac{x^\eta}{\eta}-4\frac{y^\gamma}{\gamma}}, 0<\gamma,\eta\leq 1, \tag{48}$$

with initial conditions

$$\psi\left(\frac{x^\eta}{\eta},0\right)=e^{\frac{x^\eta}{\eta}}, \psi_y\left(\frac{x^\eta}{\eta},0\right)=-2e^{\frac{x^\eta}{\eta}}, \tag{49}$$

and boundary conditions

$$\psi\left(0,\frac{y^\gamma}{\gamma}\right)=\psi_x\left(0,\frac{y^\gamma}{\gamma}\right)=e^{-2\frac{y^\gamma}{\gamma}}. \tag{50}$$

**Solution.** Applying the CDLST on (48) and the CLT to (49) and the CST to (50), to get

$$\Psi(v,\omega)=\frac{1}{(v-1)(1+2\omega)}+\frac{\omega^2}{\omega^2v^2-2\omega-1}\mathscr{L}_x^\eta\mathscr{S}_y^\gamma\left[\psi^2-e^{2\frac{x^\eta}{\eta}-4\frac{y^\gamma}{\gamma}}\right]. \tag{51}$$

Taking $(\mathscr{L}_x^\eta)^{-1}(\mathscr{S}_y^\gamma)^{-1}[\Psi(v,\omega)]$ of (51), we get

**Table 4**

The absolute error of $\psi\left(\frac{x^\eta}{\eta},\frac{y^\gamma}{\gamma}\right)$ gained by CDLSMD method for Equation (42) at different values of $\eta$, $\gamma$ and $y$.

| x | y | $\eta$ | $\gamma$ | Exact | CDLSMD | $\|\psi_{Exact}-\psi_{CDLSMD}\|$ |
|---|---|---|---|---|---|---|
| 1 | 0.1 | 0.8 | 0.8 | 0.01 | 0.0613254 | $5.13254\times 10^{-2}$ |
|   | 0.2 |   |   | 0.04 | 0.185904 | $1.45904\times 10^{-1}$ |
|   | 0.3 |   |   | 0.09 | 0.355659 | $2.65659\times 10^{-1}$ |
| 1 | 0.1 | 0.9 | 0.9 | 0.01 | 0.0241563 | $1.41563\times 10^{-2}$ |
|   | 0.2 |   |   | 0.04 | 0.084117 | $4.4117\times 10^{-2}$ |
|   | 0.3 |   |   | 0.09 | 0.174521 | $8.45212\times 10^{-2}$ |

$$\psi\left(\frac{x^{\eta}}{\eta}, \frac{y^{\gamma}}{\gamma}\right) = e^{\frac{x^{\eta}}{\eta} - 2\frac{y^{\gamma}}{\gamma}} + \left(\mathscr{L}_{x}^{\eta}\right)^{-1} \left(\mathscr{S}_{y}^{\gamma}\right)^{-1} \left[\frac{\omega^{2}}{\omega^{2} v^{2} - 2\omega - 1} \mathscr{L}_{x}^{\eta} \mathscr{S}_{y}^{\gamma} \left[\psi^{2} - e^{2\frac{x^{\eta}}{\eta} - 4\frac{y^{\gamma}}{\gamma}}\right]\right]. \tag{52}$$

Using the modified decomposition method and substituting (20) in (52) and using the findings in (22), (23), and (24), we get the components of the solution as

$$\psi_{0}\left(\frac{x^{\eta}}{\eta}, \frac{y^{\gamma}}{\gamma}\right) = e^{\frac{x^{\eta}}{\eta} - 2\frac{y^{\gamma}}{\gamma}},$$

$$\psi_{1}\left(\frac{x^{\eta}}{\eta}, \frac{y^{\gamma}}{\gamma}\right) = \left(\mathscr{L}_{x}^{\eta}\right)^{-1} \left(\mathscr{S}_{y}^{\gamma}\right)^{-1} \left[\frac{\omega^{2}}{\omega^{2} v^{2} - 2\omega - 1} \mathscr{L}_{x}^{\eta} \mathscr{S}_{y}^{\gamma} \left[A_{0} - e^{2\frac{x^{\eta}}{\eta} - 4\frac{y^{\gamma}}{\gamma}}\right]\right] = 0, A_{0} = \psi_{0}^{2},$$

$$\psi_{2}\left(\frac{x^{\eta}}{\eta}, \frac{y^{\gamma}}{\gamma}\right) = \left(\mathscr{L}_{x}^{\eta}\right)^{-1} \left(\mathscr{S}_{y}^{\gamma}\right)^{-1} \left[\frac{\omega^{2}}{\omega^{2} v^{2} - 2\omega - 1} \mathscr{L}_{x}^{\eta} \mathscr{S}_{y}^{\gamma} [A_{1}]\right] = 0, A_{1} = 2\psi_{0}\psi_{1}.$$

Following that, we get have the solution of (48) in the form

$$\psi\left(\frac{x^{\eta}}{\eta}, \frac{y^{\gamma}}{\gamma}\right) = e^{\frac{x^{\eta}}{\eta} - 2\frac{y^{\gamma}}{\gamma}}, \tag{53}$$

we can see if $\eta = \gamma = 1$; then the accurate solution is,

$$\psi(x, y) = e^{x - 2y}.$$

In the figures below, we sketch the 3D plot of the accurate solution of Equation (48) in Fig. 5 (a) which is the same approximate obtained solution, when putting $\eta = \gamma = 1$. In Fig. 5 (b), we sketch the approximate solution with various values of $\eta$ and $\gamma$: $\eta = \gamma = 1, 0.9, 0.8$.

Table 5 below, presents the absolute exact errors obtained from computing the absolute difference of the accurate and CDLSMD solutions obtained in this problem with different values of the variables $y = 0.1, 0.2, 0.3$ at $x = 1$.

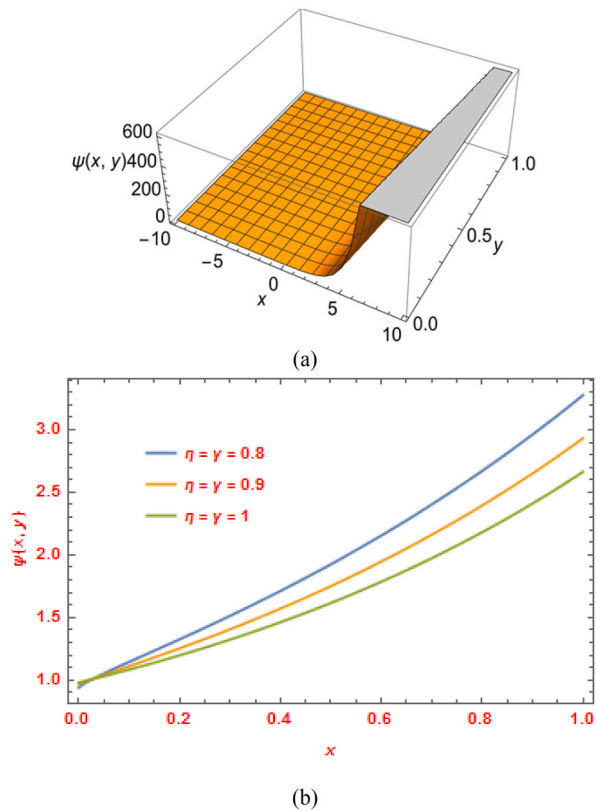

**Fig. 5.** (a) The 3D plots of the solution of Equation (53) gained by the presented method comparing to the exact solution (b) The CDLSMD solution of $\psi\left(\frac{x^{\eta}}{\eta}, \frac{y^{\gamma}}{\gamma}\right)$ for Equation (53) at $\eta = \gamma = 1, 0.9, 0.8$.

**Table 5**

The error of $\psi\left(\frac{x^\eta}{\eta}, \frac{y^\gamma}{\gamma}\right)$ g ained by CDLSMD method for Equation (53) at various values of $\eta$, $\gamma$ and $y$.

| $x$ | $y$ | $\eta$ | $\gamma$ | Exact | **CDLSMD** | $\|\psi_{Exact} - \psi_{CDLSMD}\|$ |
|---|---|---|---|---|---|---|
| 1 | 0.1 | 0.8 | 0.8 | 2.22554 | 2.3485 | $1.22959 \times 10^{-1}$ |
|   | 0.2 |     |     | 1.82212 | 1.75091 | $7.12096 \times 10^{-2}$ |
|   | 0.3 |     |     | 1.491823 | 1.34421 | $1.47616 \times 10^{-1}$ |
| 1 | 0.1 | 0.9 | 0.9 | 2.22554 | 2.29642 | $7.08757 \times 10^{-2}$ |
|   | 0.2 |     |     | 1.82212 | 1.80229 | $1.98285 \times 10^{-2}$ |
|   | 0.3 |     |     | 1.491823 | 1.43211 | $5.97162 \times 10^{-2}$ |

## 6. Results and discussion

This section compares the exact and approximate solutions utilizing tables and graphs to discuss the precision and applicability of the suggested strategy. The 3D plot solutions for Examples 1–5 generated using the current method are shown in (a) in Figs. 1–5 in contrast to the precise solutions at $\eta = \gamma = 1$. These numbers demonstrate how nearly identical the exact solutions are to the approximations obtained by CDLSMD method. Comparing the line plots of the approximate solutions from the suggested method and the exact solutions of the considered applications (1–5) for a variety values of $\eta$ and $\gamma$, is shown in (b) in Figs. 1–5. The figures show that the numerical solutions approach the exact solutions when $\eta, \gamma \to 1$. Table 1 through 5 present a comparison of each example's exact and approximate solutions in terms of absolute error at $x = 1$, and at $\eta = \gamma = 0.8$ and 0.9 for various values of $y$. As presented in the tables and figures, it is shown that the proposed method for finding the s converges swiftly towards the accurate solutions.

## 7. Conclusion

In this article, the CDLST definition has been presented. We started by using the CDLST for a few common tasks. The discussion of a few fundamental theorems and elementary properties connected to the CDLST follows. In order to demonstrate the effectiveness and usefulness of the mentioned double transform, we also used the CDLST combined with the decomposition modified approach to obtain accurate solutions of a number of classes of nonlinear differential equations in the sense of conformable derivatives. We draw the results and compare them to the accurate solutions in the integer case to demonstrate how effectively the results were attained. We conclude that the proposed approach is reliable, suitable, and efficient for obtaining precise solutions to nonlinear conformable problems. Additionally, the computations of the CDLSMD method take less time and effort and require fewer resources than those of other approaches. Thus, our goal in this study is achieved and we proved that the proposed method is applicable and accurate in handling nonlinear conformable fractional problems, and it could present exact solutions unlike other numerical methods. For possible future works, the authors will try solving nonlinear fractional integral equations.

**Funding**

This research received no external funding.

**Author contribution statement**

Shams A. Ahmed; Rania Saadeh; Ahmad Qazza; Tarig M. Elzaki: Conceived and designed the experiments; Performed the experiments; Analyzed and interpreted the data; Contributed reagents, materials, analysis tools or data; Wrote the paper.

**Data availability statement**

No data was used for the research described in the article.

**Declaration of competing interest**

The authors declare that they have no known competing financial interests or personal relationships that could have appeared to influence the work reported in this paper.

**Acknowledgments**

The authors would like to thank the editor and the kind referees, who choose to remain anonymous, for their insightful comments that helped to improve the paper's final edition.